\documentclass[12pt, leqno]{article}

\usepackage{amsmath}
\usepackage{amssymb}
\usepackage{mathrsfs}
\usepackage{stmaryrd}
\usepackage{enumerate}
\usepackage{amsthm}

 \usepackage[style=authoryear-comp,%
 				maxnames=2,%
 				mincrossrefs=2,%
 				hyperref=true,%
 				firstinits=true,
 				backend=bibtex,
 ]{biblatex}

\DeclareQuotePunctuation{}
\bibliography{button-walsh.bib}

\makeatletter
\newrobustcmd*{\parentexttrack}[1]{%
  \begingroup
  \blx@blxinit
  \blx@setsfcodes
  \blx@bibopenparen#1\blx@bibcloseparen
  \endgroup}
\AtEveryCite{%
  }
\makeatother

\usepackage[pdftex,
  bookmarks = false,
  pdfstartview = FitBH,
  linktocpage = true,
  pdfhighlight = /O,
  colorlinks=true,
  linkcolor=blue,
  citecolor=blue,
  filecolor = blue,
  urlcolor = blue,
  menucolor = blue,
]{hyperref}

\theoremstyle{definition}

\makeatletter
\newtheorem*{rep@theorem}{\rep@title}
\newcommand{\newreptheorem}[2]{%
\newenvironment{rep#1}[1]{%
 \def\rep@title{#2 \ref{##1}}%
 \begin{rep@theorem}}%
 {\end{rep@theorem}}}
\makeatother

\newcommand\model[1]{\ensuremath{\mathcal{#1}}}
\newcommand\lang[1]{\ensuremath{\mathscr{#1}}}

\newcounter{foo}

\newenvironment{listclean}
               {\list{}{\itemsep=0pt \parsep=0pt}}
               {\endlist\ignorespaces}

\numberwithin{equation}{section}
\numberwithin{enumi}{section}
\numberwithin{foo}{section}
\numberwithin{thm}{section}

\title{Structure and Categoricity: Determinacy of Reference and Truth-Value in the Philosophy of Mathematics}
\author{Tim Button\footnote{University of Cambridge, button@cantab.net} {} and Sean Walsh\footnote{Department of Logic and Philosophy of Science, University of California, Irvine, swalsh108@gmail.com or walsh108@uci.edu}}
\begin{document}
\maketitle

\begin{abstract}
\noindent This article surveys recent literature by Parsons, McGee, Shapiro and others on the significance of categoricity arguments in the philosophy of mathematics. After discussing whether categoricity arguments are sufficient to secure reference to mathematical structures up to isomorphism, we assess what exactly is achieved by recent `internal' renditions of the famous categoricity arguments for arithmetic and set theory.
\end{abstract}

\section{Introduction}\label{sec01}
In recent decades, the philosophy of mathematics  has been awash with appeals to, and reformulations of, the famous categoricity theorems of Dedekind and Zermelo. This has gone hand-in-hand with the development of various versions of structuralism. While categoricity arguments and structuralism can be  pursued independently, recent philosophical discussion has tended to link them, by considering responses to sceptical concerns about both the determinacy of reference of mathematical language, and the determinacy of truth-values of mathematical statements. 

We begin by charting considerations which have helped to threaten the determinacy of reference, and thereby pushed structure and categoricity to the forefront of philosophy of mathematics (\S\ref{sec02}). We then survey the use of Dedekind's and Zermelo's categoricity results, focussing particularly on whether they provide us with a means for grasping particular mathematical structures (\S\ref{sec03}). Finally, we turn to \emph{internal} categoricity results (\S\ref{sec04}). These have been the focus of recent literature and, by emphasising how they differ from the traditional, \emph{external}, categoricity results, we demonstrate their rather different philosophical significance.

\section{Towards structure and categoricity}\label{sec02}
We begin with one of the simplest possible observations from model theory: given a structure, it is very easy to construct an isomorphic structure. In more detail:
	\begin{listclean}
		\item[\textbf{The Push-Through Construction.}] Let~$\lang{L}$ be any signature, let~$\model{M}$ be any~$\lang{L}$-structure with underlying domain $M$, and let~$\pi:M\rightarrow N$ be any bijection. We can use~$\pi$ to induce another~$\lang{L}$-structure~$\model{N}$ with underlying domain $N$, just by `pushing through' the assignments in~$\model{M}$, i.e.\ by stipulating that~$\textbf{s}^{\model{N}} = \pi(\textbf{s}^\model{M})$ for each~$\lang{L}$-symbol~$\textbf{s}$. Having done this, one can then check that~$\pi: \model{M} \rightarrow \model{N}$ is an isomorphism.\footnote{In more detail, if $c$ is a constant symbol of $\lang{L}$, then we set $c^{\model{N}}= \pi(c^{\model{M}})$, and if $R$ is an $n$-ary relation symbol of $\lang{L}$, then we set $R^{\model{N}} = \{(\pi(a_1), \ldots, \pi(a_n)): (a_1, \ldots, a_n)\in R^{\model{M}}\}$, and if $f$ is an $n$-ary function symbol of $\lang{L}$, then we set $f^{\model{N}}(\pi(a_1), \ldots, \pi(a_n)) = \pi(f^{\model{M}}(a_1, \ldots, a_n))$.
For more discussion of this construction, see \parencite[225--31]{Button2013}.} 
	\end{listclean}\noindent
This trivial result is at the core of two considerations which push us in the direction of a focus on structure and categoricity.

\subsection{Benacerraf's use of Push-Through}\label{sec02b}
Since isomorphic copies of structures are \emph{so} easy to come by, it is no surprise that, for many mathematical purposes, it seems not to matter which of two isomorphic models one works with. 

The most famous philosophical statement of this point is due to Benacerraf. He focussed specifically on the fact that, when doing arithmetic, it makes no difference whether we think of the natural numbers as von Neumann's finite ordinals or as Zermelo's finite ordinals. 	
	\begin{quote}
		For arithmetical purposes, the properties of numbers which do not stem from the relations they bear to one another in virtue of being arranged in a progression are of no consequence whatsoever \parencite*[69-70]{Benacerraf1965}
	\end{quote}
Otherwise put: for arithmetical purposes, all we require is that the progression is both long enough and also appropriately structured. Since we can use the Push-Through Construction to induce the appropriate structure, this reduces to the requirement that we have \emph{enough} things.

Such observations suggest that mathematicians---and so, perhaps, philo\-so\-phers---can (typically) focus only on mathematical `structure'. And indeed, this attitude is reflected in some parts of mathematical practice. Mathematical discourse is rich with apparent definite descriptions of mathematical structures, like `\emph{the} natural number structure' or `\emph{the} Klein four-group'; but, given the Push-Through-Construction and related issues, it is equally  rich with the idea that we only care about the identity of such entities `up to isomorphism'. This motivates a compelling idea: \emph{mathematical structures, as discussed informally by mathematicians, are best explicated by isomorphism types}.\footnote{Unfortunately, there is an infelicity here. A model theorist's structures are always relative to a specific signature. Consequently, the natural numbers with `$0$' and `$s$' as the sole primitives induces a \emph{different} equivalence class of isomorphic structures than the natural numbers with `$0$', `$\times$' and `$+$'. There are detailed moves one could make in response to this, such as moving from \emph{isomorphism} to \emph{definitional equivalence} (See the discussion in \textcites[535--6]{Resnik1981aa}[207--8]{Resnik1997aa}[91]{Shapiro1997}[95]{Walsh2014aa}). However, we shall simply slur over this point in what follows, and treat mathematical structures (informally construed) as isomorphism types.} So, when a mathematician claims that the theory of arithmetic picks out `\emph{the} natural number structure', we shall treat this as the claim that it picks out a particular equivalence class of isomorphic models.

This is a very natural claim concerning arithmetic. But it is worth emphasising that not every mathematical theory is like the theory of arithmetic in this respect. Group theory does not aim to pick out any particular mathematical structure (informally construed): the whole \emph{point} of group theory is that it can be applied to many different mathematical areas, such as modular arithmetic, certain classes of permutations, etc. Similar points hold for theories governing rings, fields, etc. Following \textcite[40--1]{Shapiro1997}, we call these kinds of theories \emph{algebraic}. These are to be contrasted with \emph{non-algebraic} theories which, like arithmetic, aim to `describ[e] a certain definite mathematical domain' \parencite[39]{Grzegorczyk1962aa} or `specify \emph{one particular interpretation}' \parencite[273]{Kline1980aa}.

\subsection{Putnam's use of Push-Through}\label{sec02a}
We just used the Push-Through Construction to motivate a focus on mathematical structure, informally construed, which is to be formally treated in terms of isomorphism types. However, there is an equally celebrated use of the Push-Through Construction.

Philosophers sometimes use model theory to explicate the intuitive notions of truth and reference. On this approach, a model's assignment of constants to objects is taken to explicate the relation between names and their referents, and Tarski's recursive definition of satisfaction is taken to explicate the way in which statements are made true. And, with this explication as a backdrop, Putnam used model theory to raise philosophical questions about reference.\footnote{Both Putnam's permutation arguments and his just more theory manoeuvre have a long history; for details, see \textcite[p.~14 notes 1--2, p.~18 note 7, p.~27 note 2]{Button2013}.} 

Given a model which makes true everything which intuitively should be true, Putnam used the Push-Through Construction to generate a distinct but isomorphic model. (Putnam focussed specifically on the case where we generate a new model by permuting the underlying domain of the given model, rather than by substituting in new objects; hence his argument is called the \emph{permutation argument}.) Now, since the two models are isomorphic, they make true exactly the same sentences; but since they are distinct, they differ on the (explicated) reference of some symbol. And this raises the philosophical question: \emph{What, if anything, fixes reference? }

A natural reaction to Putnam's permutation argument is to claim that some models are simply more \emph{preferable} than others, as candidates for explicating reference. Howsoever `preferability' is spelled out, the crucial thought is that isomorphic models can differ as regards their preferability (see \textcites[80]{Merrill1980}[227--8]{Lewis1984aa}). But, in response, \textcites[486--7]{Putnam1977}[477]{Putnam1980}[45--8]{Putnam1981}[ix]{Putnam1983} insisted that his opponent must explain both her notion of preferability, and also why preferable structures are better at explicating reference. In doing so, Putnam maintained, his opponent would have to provide a \emph{theory} of preferability. But the permutation argument establishes that \emph{any} theory has multiple models; including this theory of preferability, coupled with anything else we want to say. So Putnam maintained that the theory of preferability is \emph{just more theory}---more grist for the permutation-mill---and cannot constrain reference. For obvious reasons, this general argumentative strategy is called the \emph{just more theory manoeuvre}.

This manoeuvre has been widely criticised (e.g.\ \textcites[225]{Lewis1984aa}[342--8]{Bays2001aa}[197--207]{Bays2008}; for further references see \textcite[p.\ 29 note 8]{Button2013}).  It is easy to see why. To use one of Putnam's own examples: suppose there are exactly as many cats as cherries, and that we are considering a permutation which sends every cat to a cherry. This permutation is supposed to make us worry that our predicate `\ldots is a cat' picks out the cherries, rather than the cats. But that permutation surely fails to respect the \emph{causal} relationships that link our use of the word `cat' to cats, rather than cherries. So, we are here suggesting that a \emph{preferable} model must respect certain causal constraints on language--object relations. In this context, the just more theory manoeuvre amounts to noting that we can reinterpret the word `causation'. But if causation does, in fact, fix reference, then this \emph{re}interpretation is just a \emph{mis}interpretation, with no philosophical significance.

In this paper, though, we are less concerned with ordinary objects, like cats and cherries, and more concerned with mathematical objects. And it is easy to see why the permutation argument and the just more theory manoeuvre might seem more threatening  in mathematical contexts (see \textcites[127, 133--5]{Hodes1984}[35--8]{McGee1997aa}). In the ordinary case, invoking causation makes it at least plausible that isomorphic models can differ as regards their preferability (though for more see \textcite[ch.\ 3--7]{Button2013}). In the mathematical case, it is much less obvious what we might invoke, if we wanted to explain why `27' refers to 27, and not to some other object.

We are not insisting, here, that these considerations provide a knock-down argument against the idea that mathematical terms determinately refer. However, together with the observations of \S\ref{sec02b}, these kinds of considerations have led some philosophers simply to embrace the idea that mathematical terms refer indeterminately.\footnote{For example: \textcite[36]{McGee1997aa} regards the question of `how mathematical terms come to have determinate referents' as `insoluble'. \textcites[84]{Balaguer1998}[73]{Balaguer1998ab} says that platonists `have to claim that while [theories like Peano Arithmetic] truly describe collections of abstract mathematical objects, they do not pick out \emph{unique} collections of such objects'. \textcite[581]{MacBride2005} similarly writes that `indeterminacy appears to be an ineliminable aspect of reference to mathematical objects'. Indeed, the prospect of referential indeterminacy is also sometimes used as evidence for the conclusion that mathematical language does not refer at all (cf.\ \textcite[139]{Hodes1984}).}

\subsection{Shapiro and Places in structures}\label{sec03b}
Indeed, the considerations of \S\ref{sec02b} might have provided us with a further reason to abandon the determinacy of reference for mathematical terms. Following \S\ref{sec02b}, suppose we maintain that our theory  of arithmetic picks out a particular isomorphism type (i.e.\ a class of isomorphic models). Then we might say that the theory of arithmetic \emph{itself} refers determinately. But when we consider a particular arithmetical term, such as `27', there is no \emph{single} object in the isomorphism type for it to pick out. So, on this view, if `27' refers at all, then it surely refers to all `the 27s' of all the isomorphic models equally, as it were. And that amounts to radical referential indeterminacy.

That said, Shapiro has prominently attempted to reconcile the focus on `structure' with the \emph{determinacy} of reference for mathematical terms. Structures, according to Shapiro---which we shall call \emph{ante-structures}---are abstract entities, consisting of \emph{places}, with certain intra-structural relations holding between them.\footnote{Shapiro's ante-structures can be compared with fruitfully with Resnik's \parencites*{Resnik1981aa}{Resnik1997aa} \emph{patterns}; his places with Resnik's \emph{positions}.} These abstract entities should be considered along the lines of universals which can be multiply realised. Shapiro calls such realisations  \emph{systems}, and a \emph{place-holder} is an object which, on that realisation, instantiates a particular place in the realised ante-structure \parencite*[73--4]{Shapiro1997}. Shapiro's \emph{ante rem} structuralism is then the claim that ante-structures exist independently from the systems realising them \parencites*[9, 84--5, 109]{Shapiro1997}. 

Now: as noted at the very end of \S\ref{sec02a}, there is no single object in the isomorphism type of the natural numbers which could serve as the preferred referent for the arithmetical term `27'. But, according to Shapiro, the \emph{ante-structure} for arithmetic contains \emph{places} which act as potential reference candidates for arithmetical terms. Thus, Shapiro maintains that `27' determinately refers to the 27-place in the ante-structure of arithmetic \parencite*[pp.~14, 55 note 15, 141 note 8]{Shapiro1997}. 

However, since places in ante-structures are abstract objects, we might well ask \emph{why} the expression `27' refers to the 27-position in the natural number structure, rather than any \emph{other} abstract entity. To make the question especially tricky, note that the Push-Through Construction will allow us to create an arithmetical structure (in the model theorist's sense, or a \emph{system} in Shapiro's terminology) according to which `27' refers to \emph{any} entity we like.  (For similar objections, see \textcites[80--4]{Balaguer1998}[193--6]{Hellman2001}[546]{Hellman2005}[151]{McGee2005aa}.) In response, Shapiro might well emphasise that the Push-Through Construction generates only a \emph{model-theoretic} structure, and not an \emph{ante}-structure, and then go on to argue that ante-structures are the \emph{preferred} reference candidates of our mathematical language. However, this would return us to the dialectic of \S\ref{sec02a}: given the abstractness of all the entities concerned, it is not immediately clear \emph{why} ante-structures should attract our words towards them.

\subsection{Further reflections}
We have considered two uses of the Push-Through Construction, both of which suggest we should focus on structure (informally construed), and (perhaps) abandon the determinacy of reference for mathematical terms. 

These two uses of the Push-Through Construction are linked technically and philosophically. However, they are also linked historically. Benacerraf's famous \cite*{Benacerraf1965}-paper grew out of his \cite*{Benacerraf1960aa}-dissertation,\footnote{The dissertation \textcite{Benacerraf1960aa} is listed as item B1 in Benacerraf's bibliography in \textcite[263]{Morton1996aa} and its relation to his papers is discussed in \textcite[24]{Benacerraf1996aa}.} in which he writes:
	\begin{quote}
		One day in conversation, Putnam made very suggestive remarks in the course of a discussion of the question `Can, or should, the numbers be identified with sets of sets?' His point was to reject the question, arguing that it arises from a distinction between number words and numbers parallel to that between, say, furniture words and furniture, and that in the former case the distinction was unwarranted. \parencite[162]{Benacerraf1960aa}
	\end{quote}
Moreover, some of their considerations were motivated by a reading of Chapter 2 of Cassirer's 1910 book, in which Cassirer defends Dedekind contra Frege on the foundations of arithmetic. Cassirer writes:
	\begin{quote}
		It is a fundamental characteristic of the ordinal theory [of natural number] that in it the individual number never means anything by itself alone, that a fixed value is only ascribed to it by its position in the total system. \parencites[47--8]{Cassirer1923aa}[62]{Cassirer1910aa}.
	\end{quote}
To our ears, this both harkens back to \textcite[\P\P{73}, 134]{Dedekind1888}, and forward to modern structuralists like Shapiro and Parsons.

However, it is worth emphasising a contrast between Benacerraf's and Putnam's considerations. Putnam forces us to ask: \emph{what does `27' refer to?} Some of Benacerraf's considerations simply ask us: \emph{which object is 27?}\footnote{Similarly, some of Benacerraf's considerations prompt us to ask: \emph{is 27 an object at all}? These are of course related, with Benacerraf writing: `I therefore argue, extending the argument that led to the conclusion that numbers could not be sets, that numbers could not be objects at all [\ldots]' \parencite*[69]{Benacerraf1965}. Obviously there are other places in Benacerraf's \cite*{Benacerraf1965}-essay where the explicit focus is on questions of reference (cf.\ \parencite*[55, 61]{Benacerraf1965}). Our point here is just that one can distill from Benacerraf's essay a particular question which is \emph{less semantic in character}.} Suppose Benacerraf-style reasons have led us to focus on structure (informally construed), i.e.\ on isomorphism types. Even if we have happily accepted that mathematical \emph{terms}, like `27', do not pick out any particular object, we will still want our mathematical \emph{theories} to pick out particular  isomorphism types. But isomorphism types are abstract objects. So we now encounter a Putnam-style concern, pitched at the level of theories rather than their terms: \emph{what could allow a theory to pick out exactly one isomorphism type?}

\section{External categoricity}\label{sec03}
The most straightforward answer to the preceding question is for the theory to be \emph{categorical}; that is, for all of its models to be isomorphic. For, if theories `pick out' their models, then categorical theories pick out a single isomorphism type.\footnote{It is worth mentioning a second reason for aiming for categoricity. Since isomorphic models satisfy exactly the same sentences, if a mathematical theory manages to be categorical, then all of its models agree on the truth-values of every sentence expressible in the theory's language.  \textcite{Corcoran1980a,Corcoran1981aa} and \textcite[18]{Awodey2002aa} note that it is not always clear whether the late-19th and early-20th century practitioners of categoricity (like Dedekind and Veblen) were motivated by categoricity \emph{per~se}, or by the distinctive kind of completeness of theories which it delivers.}

\subsection{External categoricity theorems}\label{secECT}
Unfortunately, an elementary result from model theory imposes an immediate barrier to categoricity:\newtheorem*{thmLST}{L\"{o}wenheim--Skolem Theorem}
	\begin{thmLST}
		Every countable first-order theory with an infinite model has a model of every infinite size.\footnote{For a discussion of the history of this theorem, see  \parencite[\S{4} pp. 352 ff]{Mancosu2009aa}. Proofs of this theorem are available in most textbooks on model theory, textbooks on first-order logic, and discussions of extensions of first-order logic. For instance, see: \parencite[\S{2.3} pp. 44 ff]{Marker2002}, \parencite[112]{Rautenberg2010aa}, \parencite[80]{Shapiro1991}.}
	\end{thmLST}\noindent
As an immediate corollary, no countable first-order theory with an infinite model is categorical (since isomorphic models must be of the same size). But the theories we are likely to be most interested in here---arithmetic, analysis, and set theory---all have an infinite model, and so cannot be categorical. In short, if we want categoricity, we must move beyond first-order logic. 

In this context, we might move to full second-order logic, where there are some famous categoricity results. Dedekind proved the categoricity of $\text{PA}^2$ (second-order Peano Arithmetic):
\newtheorem*{thmPAe}{Dedekind's categoricity theorem}
 	\begin{thmPAe}
      	    All full models of~$\text{PA}^2$ are isomorphic.
	\end{thmPAe}\noindent
With a little massaging, Dedekind's result also yields the categoricity of second-order real analysis (see e.g.\ \cite[84]{Shapiro1991}). A result is also available for set theory, which we discuss in \S\ref{secQuasi}. However, we shall begin by discussing the philosophical uses of Dedekind's Theorem.\footnote{For Dedekind's original proof of his theorem, see \parencites[\P{132}]{Dedekind1888}[vol. 3 p. 376]{Dedekind1930}[vol. 2 p. 821]{Ewald1996aa}. For a contemporary proof, see \parencites[82--3]{Shapiro1991}[287]{Enderton2001}.}

\subsection{Model-theoretic scepticism and just more theory}\label{secMTSJMT}
Dedekind's Theorem shows that, if we can appeal to full second-order logic, then we can axiomatise arithmetic categorically. So our question becomes: can we appeal to full second-order logic? 

The `second-order' component of `full second-order logic' is surely unobjectionable. No one can prevent mathematicians from speaking a certain way, or from formalising their theories using any symbolism they like. The qualifying expression `full', however, is rather more delicate. Here, it describes a particular semantics for second-order logic: one in which the second-order quantifiers must range over the entire powerset of the first-order domain of the structure.  But if one formalizes these semantics and their constituent notions (like powerset) in a \emph{first-order} set theory, then the L\"owenheim--Skolem result applies once again, and categoricity must be abandoned.  Equally, we can supply a \emph{Henkin} semantics for  second-order logic---essentially treating second-order logic as a form of many-sorted first-order logic---and, once again, obtain the L\"owenheim--Skolem result.\footnote{Many-sorted first-order logic is a logic in which there are different classes of variables reserved for different sorts of objects, but which is otherwise just like first-order logic. It is a natural logic in which to study e.g.\ point-line geometry, since one could have a sort reserved for points and another for lines. Many-sorted treatments of second-order logic operate with a sort for first-order objects and a sort for second-order objects, along with axiomatizations of the relation between them. For technical details, see \textcite[70--6]{Shapiro1991}, or \textcite[Chapter 6]{Manzano1996} or \textcite[\S\S{4.3}-{4.4}]{Enderton2001}.}

Such observations have been made repeatedly in the recent literature (\textcites[288]{Weston1976aa}[14--17]{Parsons1990a}[394--6]{Parsons2008}[p.\ 308 note 1]{Field1994}[319, 321--2, 338--9, 352--3]{Field2001aa}[273--5]{Shapiro2012}[120]{Vaananen2012ab}[535--42]{Meadows2013}[28]{Button2013}). Here is Putnam's statement of the general problem:
	\begin{quote}
		[\ldots] the `intended' interpretation of the second-order formalism is not fixed by the use of the formalism (the formalism itself admits so-called `Henkin models'[\ldots]), and [so] it becomes necessary to attribute to the mind special powers of `grasping second-order notions'. \parencite*[481]{Putnam1980}
	\end{quote}
The parallel to Putnam's just more theory manoeuvre is clear (see \S\ref{sec02a}). Using the L\"{o}wenheim--Skolem Theorem, we are presented with various alternatives for what our mathematical theory, $\text{PA}^2$, picks out: `the standard model' (i.e.\ one particular isomorphism type), or some `nonstandard model' (i.e.\ some other type). In trying to spell out why $\text{PA}^2$ picks out the former,  we appeal to Dedekind's Theorem. But, for the Theorem to be relevant, we must have ruled out the Henkin semantics for second-order logic; indeed, in the vocabulary of \S\ref{sec02a}, we might say that full models are \emph{preferable} to Henkin models. But the distinction between full models and Henkin models essentially invokes abstract \emph{mathematical} concepts. And these were precisely what was at issue in the first place: we were hoping to secure a grasp of such concepts by \emph{appeal to} Dedekind's Theorem. And so the worry goes that the distinction between full and Henkin models is `just more theory', and hence up for reinterpretation.

 \subsection{Model-theoretic scepticism and intermediary logics}\label{sec03c}
We are here in the presence of a kind of \emph{model-theoretic scepticism}. And so, in an attempt to mollify the sceptic, we might retreat to a logic which is weaker than full second-order logic, but still stronger than first-order logic.  One idea which has occurred to many is to:
	\begin{listclean}
       	\item[(1)] Employ a fragment of second-order logic equipped with second-order variables $X,Y,Z,\ldots$ but no second-order quantifiers, with the semantics organized so that $\varphi(X)$ comes out true precisely when it is satisfied by \emph{any} subset of the first-order domain.
		\end{listclean}
It is easy to see that we can axiomatise arithmetic categorically using such a logic: just take the normal second-order theory, $\text{PA}^2$, and delete the `$\forall X$' in front of the induction axiom.\footnote{\label{fn:SorvPrecedent}The technical point here goes back to \textcite[192--3]{Corcoran1980a}, and is also discussed by \textcite[247--8]{Shapiro1991}.  \textcites[pp. 56ff.{}]{McGee1997aa}[224--40]{Lavine1994}{Lavine1999aa}[262--93]{Parsons2008} have all invoked the framework of~(1) in defence of philosophical arguments based upon categoricity results. Importantly, though, they did so in the context of the internal categoricity results of \S\ref{sec04c}.} But that this theory is categorical is neither surprising nor very interesting (as noted by \textcites[354]{Field2001aa}[253]{Walmsley2002}[333--4]{PedersenRossberg2010}[309--10]{Shapiro2012}). For this logic is just a notational variant for the fragment of full second-order logic in which formulas begin with at most one higher-order universal quantifier and contain no further higher-order quantifiers. In this context, this notational variant is surely no more philosophically significant than the fact that, in propositional logic, we can omit the outermost pairs of brackets in a sentence without risk of ambiguity.

Other intermediate logics have been considered. The following three augmentations of first-order logic also allow for the categoricity of arithmetic:\footnote{\textcite[89--92]{Read1997} discusses options (2), (3)  and (5).  \textcites[ch.\ 9]{Field1980}[320, 338--40]{Field2001aa} defends option (3).}
	\begin{listclean}
	\item[(2)] Treat zero and successor as logical constants.
	\item[(3)] Add a new quantifier for ``there are finitely many''.
    	\item[(4)] Introduce a single one-place predicate, whose fixed interpretation in an arbitrary model of arithmetic is given by the numbers finitely far away from (the interpretation of) zero. 
	\end{listclean}
As \textcite[91]{Read1997} notes, though, invoking either (2) or (3) simply `shifts the problem from the identification of postulates characterizing [the natural numbers] categorically\ldots into the semantics and model theory of the logic used to state the postulates'. Again, we hear echoes of the just more theory manoeuvre.

A marginally more interesting approach is:
	\begin{listclean}
	\item[(5)] Add H\"{a}rtig's binary quantifier, which  expresses that there are exactly as many $\varphi$'s and $\psi$'s (for any formulas $\varphi(x)$ and $\psi(x)$).
	\end{listclean}
To produce a categorical theory of arithmetic, we now add to the usual axioms of arithmetic an axiom stating `if there are exactly as many entities less than~$x$ as there are entities less than~$y$, then~$x = y$'; such a claim would be false of certain \emph{nonstandard} numbers. But to grasp the (intended) semantics of H\"{a}rtig's quantifier seems to require a grasp of the behaviour of cardinality \emph{in general}, which again seems to presuppose within the semantics precisely the notions we were seeking to secure. Similarly, we might
	\begin{listclean}
	\item[(6)] Allow sentences containing (countably) infinitely many conjunctions and disjunctions.
	\end{listclean}
We then obtain categoricity by adding to first-order Peano arithmetic a countable disjunction saying: `everything is either zero, or the successor of zero, or the successor of that, etc.'. But to grasp this proposal, we need to grasp that `etc.'; and that looks exactly like the original challenge of grasping the (or a) natural number sequence. 

Mathematically, the most interesting alternative is:
	\begin{listclean}
	\item[(7)] Stipulate that the arithmetical function symbols `$+$' and `$\times$', even if not logical constants, must always stand for \emph{computable} functions. 
	\end{listclean}
To obtain categoricity, we then invoke Tennenbaum's Theorem \parencite*{Tennenbaum1959}, that all computable models of the first-order Peano axioms for arithmetic are isomorphic (for proofs, see \textcites{Kaye2011aa}[59]{AshKnight2000}). However, this option faces the obvious challenge of cultivating a notion of \emph{computability} that is sufficiently independent from the arithmetical notions it seeks to vouchsafe. (For more on the philosophical interpretation of Tennenbaum's Theorem, see \textcites[561--3]{McCarty1987aa}{Dean2002aa}{Dean2014aa}{Halbach2005aa}{Quinon2010aa}{ButtonSmith2012}{Horsten2012ab}.)

\subsection{Model-theoretic scepticism and transcendental arguments}\label{sec03a}
In \S\S\ref{secMTSJMT}--\ref{sec03c}, we have been dealing with a model-theoretic sceptic, who seeks to deny us anything beyond what is available in first-order logic. In response, several philosophers have recently argued that there is something deeply \emph{unstable} about such model-theoretic scepticism. There are two main lines of argument to this effect. 

The first suggests that the model-theoretic sceptic cannot even state her position without undermining it. For, the model-theoretic sceptic seems to proceed as follows: she presents us with a nonstandard model, and asks us how we can be sure that this was not the target of our mathematical language, given that it is elementary equivalent to the supposedly standard model. But, insofar as we \emph{understand} what she has presented us with, it cannot be the target of our mathematical language. To make this more concrete---although this takes us beyond arithmetic to the case of set theory---suppose the sceptic has used the L\"{o}wenheim--Skolem Theorem to present us with a merely countable model of set theory.  If that were the target of our mathematical discourse, then our word `countable' would pick out only those sets which are countable-in-the-model; but since the model itself is not countable-in-the-model, it could not now be described to us as a `countable' model. (See \textcites[287--90]{Tymoczko1989}[165--7]{Moore2001}[\S3]{Moore2011}[\S3]{ButtonBIVMT} for developments on this theme. All three authors draw inspiration from Putnam's \parencites*[487]{Putnam1977}[ch.\ 1]{Putnam1981} brain-in-vat argument.)

The second line of argument suggests that the model-theoretic sceptic is not entitled, by her own lights, to employ the first-order model theory with which she raises a sceptical challenge. For, the model-theoretic sceptic is happy to employ the full power of first-order model theory. However, as \textcite[345]{Bays2001aa} notes, `the notions of finitude and recursion are needed to describe first-order model theory', since first-order formulas `can be of arbitrary \emph{finite} length, but they cannot be infinite', and first-order satisfaction is defined recursively (cf.\ also \textcites[397, 410--11]{Field1994}[318, 338, 343]{Field2001aa} and \textcites[109]{Benacerraf1985aa}). Thus, the use of first-order model theory itself presupposes a grasp of the natural numbers (at least up to isomorphism). So if the model-theoretic sceptic employs first-order model theory in an attempt to argue that the natural numbers are not unique up to isomorphism, she saws off the branch on which she sits.

\subsection{Other attitudes towards categoricity}\label{sec03d}
Throughout \S\ref{sec03}, we have considered Dedekind's Categoricity Theorem as a means to grasp some particular isomorphism type.  However, we should remember that a categoricity theorem, as a result of pure mathematics, does not dictate its terms of philosophical employment. So it is worth surveying a variety of different attitudes that one might take, when trying to assess the potential significance of a categoricity theorem.

\emph{The Algebraic Attitude.} Imagine a character who denies that there is a single natural number structure (informally construed). Instead, she thinks of arithmetic as algebraic. Such a character need not doubt the perfect mathematical rigour of Dedekind's Theorem. Rather, she---and all parties---should conclude that any indeterminacy in the natural numbers is simply mirrored by indeterminacy in second-order quantification, so that second-order logic is inevitably algebraic (see \textcite[308]{Shapiro2012} and cf.\ \textcites[107]{Mostowski1967aa}[92]{Read1997}). 

\emph{The Infer-To-Stronger-Logic Attitude.} One might think that mathematical practice itself dictates that there is exactly one natural number structure (informally construed). If we treat this structure (informally construed) in terms of isomorphism types, then we can use the L\"{o}wenheim--Skolem Theorem to argue that we must have access to resources beyond first-order logic. Indeed, for those who are not too worried by model-theoretic scepticism, this pattern of argumentation may look like nothing more than an inference to the best explanation. This is part of Shapiro's approach to second-order logic \parencites*[pp. xii--xiv, 100, 207, 217--8]{Shapiro1991}[306]{Shapiro2012}, and Shapiro cites \textcite[p. 326 n. 535]{Church1956} as a precedent (see also \cite[36]{Benacerraf1996aa}).

\emph{The Start-With-Full-Logic Attitude.} Suppose we begin with the idea that our grasp on full second-order logic is completely unproblematic. We can then use our antecedently given grasp of this logic, coupled with our treatment of mathematical structures (informally construed) as isomorphism types, to explain how we can pin down \emph{the} natural numbers. Equally, we can explain how it is that we are equipped with a ``finitary means of characterizing the infinite'' (cf.\ \textcites[89]{Read1997}).

\emph{The Holistic Attitude.} The previous two attitudes seem to presuppose a clean separation between logic and mathematics. However, one might well think that the boundary between the two is rather artificial (see \textcites[311--22]{Shapiro2012}{Vaananen2012ab}). So the final alternative is to embrace a more holistic attitude, according to which our grasp of (e.g.\ full second-order) logic and our grasp of determinate mathematical structures (in the informal sense) come together, with each helping to illuminate the other. 

On any of the last three attitudes, the existence of a categorical axiomatization is likely to be regarded as a hallmark of \emph{success} in the project of providing axiomatisations (cf.\ \textcites[92]{Read1997}[525--7, 536--40]{Meadows2013}). Whether that view is correct or not, our point here is fairly simple. By themselves,  categoricity theorems are simply results within pure mathematics, and they are open to many different reactions. Moreover, for the theorems to have any philosophical significance, it seems that we need to have taken the step made in \S\ref{sec02b}, of identifying mathematical structures (in the informal sense) with isomorphism types (or something very close).

\subsection{The quasi-categoricity of set theory}\label{secQuasi}
We now turn from the case of arithmetic to the case of set theory. 
The crucial result here is Zermelo's proof of the \emph{quasi}-categoricity of $\text{ZFC}^2$ (second-order Zermelo--Fraenkel set theory):
 \newtheorem*{thmZe}{Zermelo's quasi-categoricity theorem}
	\begin{thmZe}
          	Given any two full models of~$\text{ZFC}^2$, either they are isomorphic, or one is isomorphic to a proper initial segment of the other.\footnote{For Zermelo's original proof, see \textcites{Zermelo1930}[pp.\ 400ff.{}]{Zermelo2010aa}[vol.\ 2 pp.\ 1219 ff.{}]{Ewald1996aa}. For a contemporary proof, see \textcite[19]{Kanamori2003aa}.}
	\end{thmZe}\noindent
Famously, \textcite[150]{Kreisel1967aa} used Zermelo's Theorem to argue that the continuum hypothesis has a determinate truth value. After all, the continuum hypothesis concerns only sets of low height, and hence is either true in all full models of $\text{ZFC}^2$, or false in all full models of $\text{ZFC}^2$, by Zermelo's Theorem.

Kreisel here seems to have adopted the Start-With-Full-Logic Attitude.  \textcite[\S5]{Isaacson2011ab} has recently followed him. But what alternative attitudes were available? 

Kreisel's remarks came in the wake of an explosion of set-theoretic independence results, initiated by Cohen's proof that the continuum hypothesis is independent from the standard axioms of set theory. And, reflecting upon these independence results, Mostowski articulated the view that:
 	\begin{quote}
         [\ldots] the incompleteness of set-theory [\ldots] is comparable [\ldots] to the incompleteness of group theory or of similar algebraic theories. These theories are incomplete because we formulated their axioms with the intention that they admit many non-isomorphic models. In [the] case of set-theory we did not have this intention but the results are just the same. \parencite[94]{Mostowski1967aa}
 	\end{quote}
In short, Mostowski argued that recent mathematical developments had shown that set theory itself was algebraic. And so, in contrast with Kreisel, Mostowski adopted a similarly Algebraic Attitude towards second-order logic; an attitude that was inevitable, given his claim that second-order logic `is a part of set theory' \parencite*[107]{Mostowski1967aa}. Hamkins has recently espoused exactly the same Algebraic Attitude: on the basis that `[s]et theory appears to have discovered an entire cosmos of set-theoretic universes' \parencite*[418]{Hamkins2012ab}, Hamkins holds that Zermelo's Theorem reveals the essentially \emph{algebraic} nature of second-order quantification \parencite*[427--8]{Hamkins2012ab}. 

We have witnessed both the Start-With-Full-Logic and the Algebraic Attitudes. But there is ammunition for the Infer-To-Stronger-Logic Attitude, even in Mostowski's own remarks. Having suggested that set theory is algebraic, Mostowski expressed a concern: 
	\begin{quote}
	     	Of course if there are a multitude of set-theories then none of them can claim the central place in mathematics. Only their common part could claim such a position; but it is debatable whether this common part will contain all the axioms needed for a reduction of mathematics to set-theory. \parencites[94--5]{Mostowski1967aa}[4]{Moore1982aa}
     \end{quote}
This suggests why independence in set theory might be thought to be philosophically troubling in a way in which, say, independence in group theory or geometry is not: it potentially threatens our contemporary image of mathematics as that which is representable within the sets.\footnote{Of course there is much more to say about Mostowski's argument; we mention it here, only because it makes vivid a reason for `worrying' about independence results.} And so one might contrapose his reasoning: If the image of set theory as foundational is thought to be sufficiently central, and genuinely is threatened by regarding set theory as algebraic, then we might be able to mount an `inference to the best explanation' that we can grasp full second-order logic. This would result in the Infer-To-Stronger-Logic Attitude. 

All that remains is the Holistic Attitude, and this is on display in Shapiro's recent work. He argues that there is a there is a `false dichotomy' between undertaking foundational research using set theory, or higher-order logic:
	\begin{quote}
		[\ldots] second-order logical consequence \emph{is} intimately bound up with set theory [\ldots]. But that does not disqualify second-order logic from logical and foundational studies. Mathematics and logic are a seamless whole, and it is impossible to draw a sharp boundary between them. \parencite[312]{Shapiro2012}
	\end{quote}
	
We have now witnessed all four attitudes towards categoricity results, in the case of set theory. However, there is one last twist to consider.  Zermelo's result is a mere \emph{quasi}-categoricity result; roughly, it says that any two full set hierarchies agree `as far as they go', but that one might outstrip the other. So, if we had hoped to explain how we grasp \emph{the} set hierarchy (i.e.\ something unique-up-to-isomorphism) by invoking the Start-With-Full-Logic Attitude together with Zermelo's Theorem, then we have a problem. 

This has been recently discussed by Isaacson.\footnote{\textcite{McGee1997aa} also discusses this limitation, but argues that there is a \emph{full} categoricity result. However, McGee (unlike Isaacson) is concerned with \emph{internal categoricity}. So we postpone discussion of his position~\S\ref{sec04e}.}  Intriguingly, from the unavailability of a complete categoricity result for set theory, Isaacson seems simply to infer that the set hierarchy \emph{itself} is somewhat indeterminate in nature. Indeed, Isaacson maintains that  `what is undecided in virtue of this degree of non-categoricity is genuinely undecided' \parencite*[p.\ 53; see also pp.\ 4, 50]{Isaacson2011ab}, and that `[t]he structure of the universe of sets is always partial and extensible' \parencite*[45]{Isaacson2011ab}. But one might ask what kind of modality Isaacson has in mind, when he uses the word `extensible'. Until this is clarified, it is natural to think that Zermelo's Theorem does not allow us to grasp \emph{the} set hierarchy, but rather, to grasp \emph{many} different set hierarchies.

\section{Internal categoricity}\label{sec04}
In recent decades, discussion of categoricity has shifted from Dedekind's and Zermelo's `classical' categoricity results, towards related \emph{internal} results. These have been at the forefront of recent work by Parsons, McGee, and Lavine. For reasons of space, we shall focus on Parsons, but we shall mention various salient contrasts to these other authors.

\subsection{Structures as predicates}\label{sec04a}

Parsons is a structuralist, and his particular brand of structuralism shapes the importance he assigns to categoricity. He defines his structuralism as:
 	\begin{quote}
     	[\ldots] the view that reference to mathematical objects is always in the context of some background structure, and that the objects involved have no more by way of a `nature' than is given by the basic relations of the structure \parencites*[40]{Parsons2008}[303, 333]{Parsons1990ab}
          \end{quote}
In short, Parsons is a \emph{contextualist} about mathematical reference, with his contexts being `background structures'. 
This allows him to respond to Benacerraf's observation (\S\ref{sec02b}) as follows. A claim like `$2=\{\{\emptyset\}\}$' has no absolute truth value, but depends upon context. So: it is true in a Zermelo-style set-theoretic context, false in a von Neumann-style set-theoretic context, and simply indeterminate in a purely number-theoretic context \parencite[p.\ 103; see also p.\ 77]{Parsons2008}.

Parsons' notion of a structure, however, takes a little longer to explain. 

When a model-theorist talks about `a model of second-order Peano arithmetic', she has in mind a triple $\langle N, 0,s \rangle$, i.e.\ a domain $N$, containing an element $0$ which implements zero, and a function $s : N \rightarrow N$ which implements the successor function. She will also consider $\text{PA}^2$, a (collection of) sentence(s) axiomatising arithmetic. Finally, she considers $\models$, which is a recursively defined, bona fide, language--object relation. And she may run this all together by saying, briefly, that $\langle N, 0, s\rangle \models \text{PA}^2$.

When Parsons talks about a `model of second-order Peano arithmetic', he has something very different in mind \parencites*[335]{Parsons1990ab}[112]{Parsons2008}. For him, this will be: a one-place predicate~$N$, a one-place function symbol~$s$, and a first-order variable~$0$, such that:\footnote{Note that this approach is possible only because~$\text{PA}^2$ is finitely axiomatisable. Of course~$\text{PA}^2$ also typically contains infinitely many instances of the Comprehension Scheme. In the present context, however, these are assumed to be relegated to the background second-order logic.}
 	\begin{listclean}
 	\item $N0 \land \forall x(Nx \rightarrow Nsx) \land \forall x(Nx \rightarrow 0 \neq sx)\land~$\\
        ~$\forall x\forall y([Nx \land Ny \land sx = sy] \rightarrow x = y) \land~$\\
        ~$\forall X([X0 \land \forall x((Xx \rightarrow Nx) \land (Xx \rightarrow Xsx))] \rightarrow \forall x[Xx \leftrightarrow Nx])$
 	\end{listclean}
 This is just a formula in pure second-order logic, which axiomatises arithmetic relative to~$N$,~$0$, and~$s$. In what follows, we abbreviate this formula by~$\text{PA}(N, 0,s)$. 

 We cannot emphasise enough the difference between the model-theorist's expression $\langle N, 0,s\rangle \models \text{PA}^2$ and Parsons' expression $\text{PA}(N, 0,s)$. Parsons' expression does not mention any (sets of) \emph{sentences}, and all suggestion of a \emph{satisfaction} relation has vanished. So, when Parsons speaks of a `model of second-order Peano arithmetic', there is no longer any hint of any language--object relation, and hence no hint of a model in the model-theorist's sense.
 
 \subsection{The internal categoricity of arithmetic}\label{sec04b}

 In common with most versions of structuralism, \textcites[13]{Parsons1990ab}[272]{Parsons2008} wants to pin down a unique natural number structure. To do this, he invokes a categoricity result: \emph{all models of second-order Peano arithmetic are isomorphic}. Stated thus, this looks almost exactly like Dedekind's Categoricity Theorem (see \S\ref{secECT}). However, we have dropped the qualifier `full', whose significance we emphasised throughout \S\ref{sec03}.  To make sense of this omission, we must recall that `models of second-order Peano arithmetic' is meant here in Parsons' sense, rather than in the model-theorist's language--object sense. So, \emph{Parsons' Categoricity Theorem} is just a sentence of pure second-order logic, of the following shape:
 	\begin{align*}
     	\forall N_1 \forall 0_1\forall s_1 \forall N_2 \forall 0_2 \forall s_2 ([&\text{PA}(N_1, 0_1, s_1) \wedge \text{PA}(N_2, 0_2, s_2)] \rightarrow {} \\
         & \exists F\, \text{Iso}(F, N_1, 0_1, s_1, N_2, 0_2, s_2))
     \end{align*}
 where `$\text{Iso}(F, \ldots)$' states that the second-order object~$F$ is a (second-order) isomorphism.\footnote{So~$\text{Iso}(F, N_1, 0_1, s_1, N_2, 0_2, s_2)$ just abbreviates the conjunction of:~$\forall x\forall y(Fx = y \rightarrow (N_1x \wedge N_2y))$, i.e.\  `$F$'s domain is~$N_1$ and codomain is~$N_2$'; and~$\forall x \forall y(Fx = Fy \rightarrow x = y) \wedge \forall y(N_2y \rightarrow \exists x\, Fx = y)$, i.e.\ `$F$ is a bijection'; and~$F0_1 = 0_2 \wedge \forall x (N_1x \rightarrow Fs_1 x = s_2 Fx)$, i.e.\ `$F$ preserves structure'.} At the risk of repetition then: whilst we might offer similar informal glosses on Parsons' and Dedekind's Theorems, they are very different statements. Dedekind's Categoricity Theorem concerns structures, satisfaction and sentences; Parsons' Theorem is just a sentence of pure second-order logic. 

 Crucially, Parsons' Categoricity Theorem can be proved within a \emph{deductive} system for second-order logic (for a simple proof, see \cite{Vaananen2015aa}). Consequently, we do not need to say anything about the \emph{semantics} for second-order quantification in order to secure Parsons' result.

We call Parsons' result an \emph{internal} categoricity result, and contrast this with Dedekind's \emph{external} categoricity result.\footnote{\textcite[249--51]{Walmsley2002} seems to be the first author to use the phrase `internal categoricity', but \textcite[31--9]{Parsons1990a} seems to be the first author to invoke such a result. Sometimes internal categoricity results are called \emph{relative} categoricity results.}  These labels are appropriately suggestive. Crudely: an external result requires that we stand back from the object language we were using, and consider its semantics in some metalanguage. By contrast, an internal result can be proved \emph{within} the object language.

 \subsection{What might an internal categoricity result show?}\label{sec04c}
Internal categoricity results provide us with a notion of categoricity which does not require any discussion of the semantics for full second-order logic. Given what was said in \S\ref{sec03}, this might be thought a great advantage. However, we must go slowly when considering what such results might show. 

To be very clear: internal categoricity results cannot possibly aid in the project of explaining how the semantics for mathematical language(s) get fixed, or in defeating the model-theoretic sceptic. Parsons' Categoricity Theorem is a sentence of pure second-order logic, which does not mention any sentences, and so tells us nothing (directly) about any language--object relation, or indeed about \emph{any} semantic facts. Of course, we might stand back from Parsons' Theorem, as it were, and start to explore the semantics for the language in which we proved that Theorem. In so doing, we will see that Parsons' Theorem entails Dedekind's Theorem \parencites[64]{Lavine1999aa}[99]{Vaananen2012ab}. But now Parsons' Theorem has simply taken us on a long-winded detour back to \S\ref{sec03}.

So let us now consider how Parsons himself attempts to use internal categoricity. \textcites[35--8]{Parsons1990a}[283--8]{Parsons2008} imagines two characters, Kurt and Michael, who are both `doing arithmetic'. Kurt has a predicate `\ldots is a natural number', which we can symbolise as~$N_{\emph{kurt}}$ or~$N_k$, he has a name `zero', which we can symbolise as~$0_k$, and he has a notion of function `the successor of\ldots', which we can symbolise as~$s_k$. We similarly symbolise Michael's arithmetical vocabulary with~$N_m$,~$0_m$ and~$s_m$. Assuming reasonable levels of communication between Kurt and Michael, it will be obvious to them that they are engaged in \emph{somewhat} similar practices. But what they might want to show is that their languages are essentially identical; that, for arithmetical purposes, they differ only in the subscripts we have imposed; that they are `syntactically isomorphic', to use Lavine's \parencite*[47]{Lavine1999aa} phrase. 

Parsons suggests that Kurt and Michael can establish this as follows. Let us allow that Kurt and Michael are in communication with one another to the point that both are able to take the other's vocabulary into his own language (and both know this). Both can now prove Parsons' Categoricity Theorem; and so, since they both have access to each other's vocabulary, both can prove:\footnote{\label{fn:WhySorv}In footnote \ref{fn:SorvPrecedent}, we noted that Parsons, Lavine and McGee all invoke the first weaker logic~(1) in our enumeration of such logics from \S\ref{sec03c}. We can now see why this is so. Parsons, Lavine and McGee want the predicate-variables of their logic to be open ended over any future expansion of the object language. This is why Kurt is able to establish a result which uses a predicate which originally belonged to Michael's vocabulary (and vice versa). However, see \textcites[358--60]{Field2001aa} for discussion of Kurt's warrant for assuming that Michael may countenance and perform induction on predicates which Kurt introduces.}
 	\begin{align*}
     	[\text{PA}(N_k, 0_k, s_k) \wedge \text{PA}(N_m, 0_m, s_m)] \rightarrow \exists F\, \text{Iso}(F, N_k, 0_k, s_k, N_m, 0_m, s_m)
     \end{align*}
Furthermore, both can presumably see that the antecedent obtains: they affirm one of the conjuncts themselves, and their interlocutor happily affirms the other. They therefore obtain the consequent. And this guarantees that, for arithmetical purposes, their languages differ only in the subscripts we have imposed. 

 To be clear, the foregoing argument will not decide whether~$0_k = 0_m$. This is no surprise, given what we said at the start of \S\ref{sec04a}. Indeed, context may well leave it indeterminate whether `$0_k = 0_m$' is true. But this does seem to generate a slight tension for Parsons, since a natural kind of disquotation principle will tell us that it may be indeterminate whether~$0_k = 0_m$; and yet it is a theorem that~$0_k = 0_m \vee 0_k \neq 0_m$. The more general point here is as follows: Parsons' notion of `structure' eschews all semantic notions; however, his contextualism is essentially and deliberately semantic; and we are not given much guidance as to how the non-semantic and the semantic relate to one another.
 
 We can illustrate the same point by considering the impact of Parsons' Theorem on the claim that all sentences of pure arithmetic have determinate truth-values. Once Kurt and Michael have established the existence of their second-order isomorphism, they can see that if they ever disagree (modulo subscripts) about any arithmetical sentence, then only one of them is right. But does this enable them to say that every sentence in the language of arithmetic is either true or false? Presumably we shall want to say so;\footnote{In \S\ref{sec04e}, we shall see that McGee appeals to the internal categoricity of set theory in an attempt to show that every sentence in the language of (pure) set theory has a determinate truth value \parencite*[40--2]{McGee1997aa}.} but this evidently requires the use of semantic principles which do not themselves appear anywhere within the categoricity result.

If the semantic principles that we need to invoke here are, in the end, just the notions provided by model theory, then Parsons' novel understanding of structures, and his invocation of internal categoricity, would have been a needless detour back to \S\ref{sec03}. Fortunately, Parsons seems to have a rather different view in mind. Roughly put, he holds that, once we consider semantic questions \emph{for} our mother tongue, we can only answer them by `acquiescence in our mother tongue' \parencites*[39]{Parsons1990a}[288]{Parsons2008}. This seems to amount to the suggestion that there is no need to look to model theory to supply a semantic theory for the mother tongue. Rather, such a `semantic theory' merely amounts to using my own words to say what my own words refer to; at which point, all I can provide is disquotational platitudes.\footnote{Continuing with the preceding footnote: McGee does not discuss `acquiescence in our mother tongue'; however, he is a fairly thorough-going disquotationalist, and so is likely to reach a similar conclusion.} 

To investigate this notion of `acquiescence' in any detail would take us very far afield; so we shall simply grant its coherence. To close our discussion of Parsons, we shall explain how Parsons-style structuralism deals with the challenge of referential indeterminacy and model-theoretic scepticism. 

Given his contextualism, Parsons will not expect context-insensitive referential determinacy for arithmetical terms: `27' may refer to one set in one context, and to a different set in another context. So, if the challenge of referential determinacy is to pose any problems, they will have to concern referential principles which plausibly \emph{should} be context-insensitive, e.g.\ that `27' refers to 27 (or maybe that `27' refers to 27$_k$). But for Parsons, this will be guaranteed just by acquiescing in the mother tongue. 

Exactly the same is true of Parsons' response to model-theoretic scepticism. His response does not turn upon his distinctive notion of structure, or his internal categoricity result. Instead, his answer to that challenge involves acquiescing in the mother tongue. His notion of structure, and his categoricity result, come in to play further downstream, in showing that different parties are engaged in the same practice. And note that {Kurt} cannot demonstrate the existence of a \emph{shared} practice just by acquiescing in \emph{his} mother tongue; for what he needs to show is that he \emph{shares} a tongue with Michael (in the relevant sense). 

The idea is, therefore, as follows. \emph{I} can deal with the  sceptic by acquiescing in my mother tongue. But now the threat of mathematical solipsism arises. This latter threat can be addressed by invoking internal categoricity. \emph{We} thereby arrive at the possibility of mathematical \emph{intersubjectivity}. However, this falls short of establishing mathematical \emph{objectivity}, in the sense that nothing so far establishes that we are all talking about the same objects. This is ultimately how we understand Parsons' \parencites*[38--9]{Parsons1990a}[287--8]{Parsons2008} invocation of internal categoricity.\footnote{We suggest a fairly similar reading of \textcite[89]{Pollard2007}. It is not clear whether this is quite the right reading of McGee's \parencites*{McGee1997aa} invocation of internal categoricity (more on McGee below).}

\subsection{The internal categoricity of set theory}\label{sec04e}
We have focussed on internal categoricity for arithmetic. Unsurprisingly, there is an internal quasi-categoricity result for set theory (cf.\ \textcites[55--66, 89--95]{Lavine1999aa}[\S4]{Vaananen2015aa}). But of course, since it is an \emph{internal} result, it simply abbreviates a second-order formula of roughly the following shape:\footnote{So `$\text{ZFC}(V, E)$' abbreviates a fairly long conjunction, axiomatizing second-order set theory relative to `$V$' (intuitively, the sets) and `$E$' (intuitively, membership). Likewise `$\text{Seg}(F, \ldots)$' indicates that~$F$ is a (second-order) isomorphism between an `initial segment' of one and the `whole' of the other.}
 	\begin{align*}
	     	\forall V_1 \forall E_1 \forall V_2 \forall E_2([&\text{ZFC}(V_1, E_1) \wedge \text{ZFC}(V_2, E_2)] \rightarrow {}\\
		&\exists F\, \text{Seg}(F, V_1, E_1, V_2, E_2))
     \end{align*}
Crucially, and as in the case of Parsons' Categoricity Theorem, no satisfaction relation appears anywhere in this sentence. Consequently, it tells us nothing directly about any language--object relation. Nonetheless, someone might, Parsons-style, invoke internal quasi-categoricity in order to generate the intersubjectivity of set theory, `so far as it goes'. 

This squares well with Martin's aims for quasi-categoricity. On Martin's view, sets derive ultimately from properties through a transfinite process of set-formation \parencites*[112]{Martin1970aa}[6]{Martin2001aa}. But this leaves space for a Parsons-style worry: whereas Kurt forms a given property into \emph{this} set, might Michael form the same property into a \emph{distinct} set? To address this worry we might then invoke internal quasi-categoricity, Parsons-style, in order to convince us that there is no more to set-theory (at low levels) than what is \emph{common} to Kurt's and Michael's  processes (cf.\ \parencites*{Martin2001aa}[365~ff]{Martin2005aa}). 

In fact, though, we may be able to go beyond mere \emph{quasi}-categoricity. Typical formulations of set theory only involve quantification over sets. But if set theory is to be applicable, then we shall want to be able to form sets from entities that are not themselves sets, also known as \emph{urelements} (cf.\ \textcites[49]{McGee1997aa}[vi, 24, 50--1]{Potter2004}). We can then distinguish \emph{pure} sets---those which, intuitively, involve no urelements anywhere in their construction---from impure sets. Let~$\text{ZFCU}^2$ be second-order Zermelo--Fraenkel set theory with Choice and with urelements, coupled with the axiom that there is a set containing all the urelements. \textcite{McGee1997aa} has proved the internal categoricity result: \emph{any two models of ~$\text{ZFCU}^2$ with unrestricted first-order quantifiers have isomorphic pure sets.} At the risk of pedantry, since this is an \emph{internal} categoricity result, it again has the shape:\footnote{Here `$\text{ZFCU}(U, E)$' abbreviates a fairly long conjunction axiomatizing second-order set theory with urelements relative to `$U$' (intuitively, the urelements) and `$E$' (intuitively, membership), with \emph{no} explicit relativisation on the hierarchy itself (as it were). Likewise `$\text{PureIso}(F, \ldots)$' indicates that there is a (second-order) isomorphism between `the pure sets'.}
 	\begin{align*}
     	\forall U_1 \forall E_1 \forall U_2 \forall E_2([&\text{ZFCU}(U_1, E_1) \wedge \text{ZFCU}(U_2, E_2)] \rightarrow {}\\
	& \exists F\, \text{PureIso}(F, U_1, E_1, U_2, E_2))
     \end{align*}
So now we must investigate the significance of McGee's Theorem. 

McGee claims that his Theorem shows that every sentence in the language of set theory has a determinate truth value. But this is not immediately clear. Since McGee's Theorem is an \emph{internal} categoricity result, it mentions no sentences, and involves no satisfaction relation.\footnote{\textcite[49--50]{McGee1997aa} is explicit on this point. He is pushed in this direction by Tarski's Indefinibility Theorem (roughly, that no sufficiently rich, consistent theory can define a truth predicate which provably has the properties we would want it to have). As a consequence, he faces a choice: either offer a \emph{third-order} definition of satisfaction suitable for \emph{second-order} theories, or relocate to internal categoricity. Since the former approach seems to embark on a futile regress, McGee embraces the latter approach, describing it as `devious' \parencite*[47--8, 50, 62]{McGee1997aa}.} So how could it tell us about the \emph{truth} of sentences? 

Our explanation returns, once again, to our discussion of Parsons. Imagine Kurt and Michael both `doing $\text{ZFCU}^2$'. Subscripting in the obvious way, McGee's Theorem proves that $\phi_k \leftrightarrow \phi_m$, for any sentence $\phi$ in the language of $\text{ZFCU}^2$ (restricted to pure sets). Otherwise put: if set theorists ever disagree (modulo subscripts) about any (purely) set-theoretic sentence, then only one of them is right. And this is just what McGee is pointing to, when he talks about the `determinacy of truth value'.

Now, since McGee's Theorem is an internal categoricity result, we can prove the result with a suitable deductive system for second-order logic, ignoring all semantic questions. So, just as in the case of Parsons: the result does not turn upon whether we are entitled to \emph{full} second-order quantifiers; we only need access to the normal deductive system for second-order logic. Equally, the proof does not turn on whether a first-order \emph{domain} can be `absolutely unrestricted', or whether we must believe that domains are `indefinitely extensible'; we need only \emph{use} the logic, without imposing any restriction on our first-order quantifiers within our object language. However, it is worth pausing on both of these points.

First: consider again the use of second-order logic. We do not need to discuss its semantics. But, in the context of the iterative conception of set, one might simply object to the use of second-order logic \emph{even in} the object language. As Reinhardt once put it: using second-order logic in this context rather makes it seems like one just ``forgot'' to add a level of sets (\parencite*[196]{Reinhardt1974aa}, cf.\ \textcite[546]{Burgess1985aa}).

Second: consider the issue of restricted first-order quantification. When Kurt aims to show that his language is `syntactically isomorphic' to Michael's, he must establish that Michael is speaking so that $\text{ZFCU}(U_m, E_m)$. And in order to do this, Kurt must not regard Michael's first-order quantifiers as tacitly restricted, relative to Kurt's own, when he incorporates Michael's vocabulary into his own. But there is room to doubt whether Kurt can ever be warranted in treating Michael's first-order quantifiers as unrestricted in this way. (This feeds into Incurvati's \cite*{Incurvati2015aa} critique of McGee's Theorem.)\footnote{Incurvati phrases his objection in terms of domains; though see \parencite*[note 8]{Incurvati2015aa}.}

The significance of this latter objection may, however, depend upon the standard of warrant in play. We have insisted that McGee's Theorem, like Parsons', cannot help to deal with model-theoretic scepticism, since it is an \emph{internal} categoricity result. And, if we are no longer concerned to rebut a \emph{sceptic}, then Kurt may simply be able to take Michael, at his homophonic word, that he is talking about `everything, without restriction' (cf.\ \cite{Pollard2007}).

 \section{Conclusion}\label{sec05}
 
We began our survey with the traditional concerns of Putnam and Benacerraf about determinacy of reference in the case of mathematical language. To assuage this concern, a natural move is to suggest that reference in the case of mathematical language is secured only up to isomorphism. And the traditional mechanism for securing such reference is a categoricity theorem, which says that any two models of a theory are isomorphic. 

Given the L\"owenheim-Skolem theorem, any categoricity theorem for theories with infinite models has to be done in an enrichment of first-order logic. But the appeal to such enriched logics has been thought to be question-begging, since the notions that they involve, such as quantification over all subsets of the first-order domain of the model, seem just as mathematical in character as the notions which a categoricity result sought secure in the first place.

This concern, among others, motivates a transition to internal categoricity results, in which one seeks to operate in the object-language of a background second-logic and prove results that have the consequence that any two models agree on the truth-values of any sentence expressible in their common language. But it is not clear whether these approaches can successfully explicate semantic notions---such as truth and reference---using only object-language resources. For this reason, some authors suggest that internal categoricity results are valuable, not because they guarantee the determinacy of truth-values, but because they guarantee mathematical intersubjectivity.

\section{Acknowledgements}

This paper arose from a seminar series on philosophy and model theory that we ran in Birkbeck in 2011. We have since presented the paper several times: in June 2014 at PhilMath Intersem 5 at Universit\'e Paris Diderot -- Paris 7; in November 2014 at the Southern California History and Philosophy of Logic and Mathematics Group; and in November 2014 at the Logic Seminar at Cambridge University. It has been measurably improved by the generous comments and feedback which we received from the editors and referees of this journal, as well as from: Sarah Acton, Andrew Arana, Bahram Assadian, John Baldwin, Neil Barton, Timothy Bays, Liam Bright, Adam Caulton, John Corcoran, Radin Dardashti, Walter Dean, William Demopoulos, Michael Detlefsen, Sean Ebels-Duggan, Fiona Doherty, Salvatore Florio, J. Ethan Galebach, Michael Gabbay, Peter Gibson, Owen Griffiths, Bob Hale, Emmylou Haffner, Jeremy Heis, Simon Hewitt, Kate Hodesdon, Will Hendy, Wilfrid Hodges, Luca Incurvati, Eleanor Knox, {\O}ystein Linnebo, Kate Manion, Pen Maddy, Toby Meadows, Jonathan Nassim, Fredrik Nyseth, Jonathan Payne, Sara Parhizgari, Michael Potter, Paula Quinon, Erich Reck, Sam Roberts, Marcus Rossberg, Gila Sher, J\"onne Speck, Sebastian Speitel, Stewart Shapiro, Robert Trueman and John Wigglesworth.

\printbibliography

\newpage

\end{document}